\documentclass[10pt]{article}
\usepackage{graphicx,amsmath,amsfonts,amssymb,amscd,bezier,amstext}
\usepackage{fancyheadings}
\usepackage{amsthm}
\usepackage{bm}
\usepackage{color}
\usepackage{fancybox}
\usepackage{latexsym}
\usepackage{graphicx}
\usepackage[all]{xy}
\usepackage{longtable}

\begin{document}

\newtheorem{definition}{Definition}[section]
\newtheorem{theorem}{Theorem}[section]
\newtheorem{corollary}{Corollary}[section]
\newtheorem{lemma}{Lemma}[section]
\newtheorem{proposition}{Proposition}[section]
\newtheorem{example}{Example}[section]
\newtheorem{remark}{Remark}[section]

\newfont{\bms}{msbm10 scaled 1200}
\def\matR{\mbox{\bms R}}
\def\matF{\mbox{\bms F}}
\def\matZ{\mbox{\bms Z}}
\def\matN{\mbox{\bms N}}
\def\matQ{\mbox{\bms Q}}
\def\matK{\mbox{\bms K}}
\def\matL{\mbox{\bms L}}
\def\matC{\mbox{\bms C}}
\def\matM{\mbox{\bms M}}
\newcommand{\B}{b\!\!\!b}
\newcommand{\A}{e\!\!\!a}

\title{Properties of $\cal{G}$-Equivalence of Matrices}

\author{{ Miriam da Silva Pereira } \\
{\small $$}{\footnotesize Centro de Ci\^{e}ncias Exatas e da Natureza- UFPB}\\
{\small $$}{\footnotesize Cidade Universit\'aria}\\
{\small $$}{\footnotesize e-mail: miriam@mat.ufpb.br - Jo\~{a}o
Pessoa - PB Brazil}}

\maketitle

\begin{abstract}

The theorem of Hilbert- Burch provides a description of
codimension two determinantal varieties and their deformations in
terms of their presentation matrices. In this work we use this
correspondence to study properties of determinantal varieties, based
on methods of singularity theory. We
establish the theory of singularities for $n\times p$ matrices
extending previous results of Bruce and  Tari \cite{Farid} 
and Fr\"uhbis-Kr\"uger \cite{FK}. The main
result of this work is the description of equivalent conditions to
${\cal{G}}$-finite determinacy of the presentation matrix of
Cohen-Macaulay varieties of codimension $2$. We apply the results to
obtain sufficient conditions for topological triviality of deformations
of weighted homogeneous matrices.

\end{abstract}

\section{Introduction}

Properties of determinantal varieties have been studied both in
geometry and commutative algebra (see \cite{Bruns}, \cite{Bruns1}).
If $M$ is a matrix with entries in $R$ and $t=\min\{n,p\}$, J. Eagon
proved in \cite{Ea1} that the ideal $I$ generated by the $t\times t$
minors  of $M$ is Cohen-Macaulay.

In particular, in the case where $X$ is a codimension two
determinantal variety, we can use the Hilbert-Burch's theorem to
obtain a good description of $X$ and its deformations in terms of
its presentation matrix. In fact, if $X$ is a codimension two Cohen-
Macaulay variety, then $X$ can be defined by the maximal minors of a
$n\times (n+1)$ matrix. Moreover, any perturbation of a  
$n\times (n+1)$ matrix gives rise to a deformation of $X$ and any
deformation of $X$ can be obtained through a perturbation of the
presentation matrix.

In this work, we use this correspondence to study properties of
codimension two  Cohen-Macaulay varieties through their presentation
matrix. We define a group ${\cal{G}}$ acting in the space of the
matrices and use this equivalence to classify their singularities.

In \cite{FK}, A. Fr\"uhbis-Kr\"uger considers the case of matrices
defining codimension two Cohen-Macaulay singularities to study
space curves. In this case, it is possible to get the normal module
and the space $T^1$ of the deformations of first order in terms of
presentation matrices. The Tjurina number, $\tau(X)$, of the
determinantal variety $X$ is the codimension of the space $T^1$ and
coincides with the ${\cal{G}}_e$- codimension of $M$. In \cite{FKn},
A. Fr\"uhbis-Kr\"uger and A. Neumer obtain a complete list of simple
Cohen-Macaulay codimension $2$ singularities.

Quadratic matrices have been previously studied by V. I. Arnold in
\cite{arn}. Moreover, using the ${\cal{G}}$-equivalence of matrices,
J. W. Bruce, F. Tari and G. J. Haslinger obtain classifications of
simple germs of families of symmetric, skew-symmetric and square
matrices using the $\cal{G}$-equivalence (\cite{Bruce},
\cite{Farid}, \cite{has}).

In the first part of the paper, we develop the theory of
singularities for $n\times p$ matrices, extending the results of
\cite{FK}, \cite{Bruce} and \cite{Farid} related to the
infinitesimal and geometric characterization of finite determinacy
and the theorem of versal unfolding.

We give equivalent of conditions to  ${\cal{G}}$-finite determinacy
of the presentation matrix of codimension $2$ Cohen-Macaulay
varieties with isolated singularity. The corresponding result for
the contact group ${\cal{K}}$ was obtained by T. Gaffney in
\cite{Gaft}. As application of these results we study the
${\cal{G}}$-topological triviality of families of matrices.

\section{Notations and Basic Definitions}

Let $Mat_{(n,p)}(\matC)$ be the set of all $n\times p$ matrices with
complex entries, $\Delta_t\subset Mat_{(n,p)}(\matC)$ the subset
formed by matrices that have rank less than $t$, with $1\leq t\leq
\min(n,p)$. It is possible to show that $\Delta_t$ is an irreducible
singular algebraic variety,   of codimension $(n-t+1)(p-t+1)$ (see
\cite{Bruns}). Moreover the singular set of $\Delta_t$ is exactly
$\Delta_{t-1}$. The set $\Delta_t$ is called \textit{generic
determinantal variety}.

\begin{definition}
Let $M=(m_{ij}(x))$ be a $n\times p$ matrix whose entries are
complex analytic functions on $U\subset\matC^r$, $0\in U$ and $f=(f_1,\ldots\,f_q)$
defined by the $t\times t$ minors of $M$. We say that
$X$ is a determinantal variety if $X$ is defined by the equation $f_1=\ldots=f_q=0$ and codimension of $X$ i equal to $(n-t+1)(p-t+1)$.
\end{definition}

We can look at a matrix $M=(m_{ij}(x))$ as a map
$M:\matC^r\longrightarrow Mat_{(n,p)}(\matC)$, with $M(0)=0$. Then
the determinantal variety in $\matC^r$ is the set
$X=M^{-1}(\Delta_t)$ when codimension of $X$ is equal to $(n-t+1)(p-t+1)$.

Let ${\cal{O}}_r$ be the ring of  germs of analytic functions on
$\matC^r$ and $\mathcal{M}$ its maximal ideal. We denote by
$Mat_{(n,p)}({\cal{O}}_r)$ the set of all matrices $n\times p$ with
entries in ${\cal{O}}_r$. This set can be identified with ${\cal{O}}^{np}_ r$ the free module of
rank $np$.

We concentrate our attention in this paper to codimension $2$
determinantal singularities and their deformations. The following
proposition follows from the Auslander- Buschsbaum formula and the
Hilbert-Burch's Theorem.

\begin{proposition}(\cite{FK}, pg. 3994)
\begin{itemize}
\item [1) ]Let $M$ be a matrix $(n+1)\times n$ with entries in ${\cal{O}}_r$
and $f=(f_1,\ldots, f_{n+1})$ its maximal minors and, by abuse of
notation, the ideal generated by them. If $codim(V(f))\geq2$ the
following sequence
$$0\longrightarrow ({\cal{O}}_r)^n\longrightarrow
({\cal{O}}_r)^{n+1}\longrightarrow {\cal{O}}_r/(f)\longrightarrow
0$$ is exact. Moreover, ${\cal{O}}_r/(f)$ is Cohen-Macaulay and
$codim(V(f))=2$.
\item [2) ] If $X\subset\matC^r$ is Cohen-Macaulay, $codim(X)=2$ and
$X=V(I)$, then ${\cal{O}}_r/I$ has a minimal resolution of the type
$$0\longrightarrow {\cal{O}}_r^n\longrightarrow
({\cal{O}}_r)^{n+1}\longrightarrow {\cal{O}}_r/I\longrightarrow 0.$$
Moreover, there is an unit $u\in{\cal{O}}_r$ such that $I=u\cdot
f$, where $f$ is again the ideal of the  maximal minors of $M$.
\item[ 3.] Any deformation of $M$ is a deformation of
$X$;
\item[4. ] Any deformation of $X$ can be generated by a
perturbation of the matrix $M$.
\end{itemize}
\end{proposition}

It follows from this proposition that any deformation of a
codimension $2$ Cohen-Macaulay variety may be given as a
perturbation of the presentation matrix. Therefore, we can study
these varieties and their deformations using their presentation
matrices.

Given a matrix $M\in Mat_{(n,p)}({\cal{O}}_r)$, let $C_{ij}(M)$
(respectively $R_{lk}(M)$) be the matrix that has the $i$-th column
(respectively the $l$-th row) equal to the $j$-th column of $M$
(respectively the $k$-th row) and zeros in any other position.

We denote by $J(M)$ is the submodule generated by the matrices of the form $\dfrac{\partial M}{\partial x_i}$ for $1\leq i\leq r$ and ${\cal{H}}=\textrm{GL}_{ \textit{p} } ({\cal{O}}_r)\times
\textrm{GL}_{\textit{n}}({\cal{O}}_{r})$, where $\textrm{GL}_{ \textit{j} } ({\cal{O}}_r)$ denotes the group of invertible $j\times j$ matrices with entries in ${\cal{O}}_{r}$.
Let $\cal{R}$ be the group of changes of coordinates in
$(\matC^r,0)$, that is, $\cal{R}$  is the group of analytic diffeomorphism germs in $(\matC^r,0)$.

Let $M$ be a $n\times (n+1)$ matrix with entries in the maximal
ideal of ${\cal{O}}_r$, $n> 1$ and $r>1$. We denote by
$f=(f_1,\ldots,f_{n+1})$ the ideal generated by the $n\times n$
minors  of the matrix $M$ and $X$ the variety defined by $f$. The index of each $f_i$ indicates the
column removed from $M$ to compute the minor. 
As $ X $ is not in general a complete intersection, there are relations between the
component functions of $f$. Thus, the maximal rank of the jacobian matrix of $f$ is
given by $d\leq\min\{n+1,r\}$.

\section{Singularity Theory of Matrices}

Given two germs of singularities of matrices, we are interested in
studying these germs according to the following equivalence relation.

\begin{definition}
Let ${\cal{G}}={\cal{R}}\ltimes {\mathcal{H}}$ the semi-direct product of ${\cal{R}}$ and ${\mathcal{H}}$. We say that two germs $M_1,
\,\ M_2\in Mat_{\textit{(n,p)}}({\cal{O}}_r)$ are
$\cal{G}$-equivalent if there exist $(\phi,R,L)\in
{\cal{G}}$ such that $M_1=L^{-1}(\phi^*M_2)R$.
\end{definition}

It is not difficult to see that $\cal{G}$ is one of Damon's
geometric subgroups of $\cal{K}$, hence a consequence of Damon's result (\cite{Damon})
we can use the techniques of singularity theory, for instance, those
concerning finite determinacy. The notions of ${\cal
{G}}$-equivalence and $ {\cal {K}} _ {\Delta} $-equivalence, where $
\Delta $ consists of the subvariety of matrices of rank less than
the maximal rank \cite {Damon}, coincide for finitely determined
germs (see \cite{Bruce}).

The next proposition characterizes the ${\cal{G}}$-tangent space of
a matrix $M$. The proof is analogous to the proof in
\cite{Bruce} in the case of symmetric matrices.

\begin{proposition}
  \begin{enumerate}
    \item The $\cal{R}$-tangent space to the orbit of an element of
     $M\in\textit{Mat}_\textit{(n,p)}({\cal{O}}_r)$ is the ${\cal{O}}_r$-module generated by $x_j\dfrac{\partial M}{\partial x_i}$, $1\leq i,\,j\leq r$.

    \item The tangent space to the orbit of $M$ under the action of the subgroup
    ${\cal{H}}=\textrm{GL}_{ \textit{p} } ({\cal{O}}_{r})\times
\textrm{GL}_{\textit{n}}({\cal{O}}_{r})$ is the ${\cal{O}}_r$-module
generated by $C_{ij}(M)$, $1\leq i,\,j\leq p$, and $R_{lk}(M)$, $1\leq l,\,k\leq n$. 
  \end{enumerate}
\end{proposition}
It follows from the previous discussion that the $\cal{G}$-tangent
space to a germ $M$ is given by {\begin{equation*}
{\cal{T}}{\cal{G}}M={\mathcal{M}}J(M)+{{\cal{O}}_r}\{R_{lk}(M),C_{ij}(M)\},
\end{equation*}}
where $\,1\leq l,k\leq n$ and $1\leq i,j\leq p$.

Given $M\in Mat_{(n,p)}({\cal{O}}_r)$ we consider the map 
\begin{equation*}
\begin{array}{r} g_M=g:Mat_{(p,p)}({\cal{O}}_r)\times
Mat_{(n,n)}({\cal{O}}_r)\longrightarrow Mat_{(n,p)}({\cal{O}}_r)
\end{array}
\end{equation*}
given by $g(A,B)= BM+MA$.

Then it is possible to rewrite the expression of the tangent space as
\begin{equation*}
{\cal{T}}{\cal{G}}M={\mathcal{M}}J(M)+{\cal{O}}_r\,Im(g).
\end{equation*}
This equivalence relation is useful to classify determinantal singularities
and to study their deformations.

The next propositions express the normal module ${\mathcal{N}}_X$
and the space of the first order deformations $T_{X}^1$, in terms of
matrices, hence we can treat the base of the semi-universal
deformation using matrix representation (see \cite{FKT} for the
definitions of ${\mathcal{N}}_X$ and $T_{X}^1$).

\begin{proposition}(\cite{FK}, pg. 3996)
Let $M$ be a $(n+1)\times n$ matrix with entries in the maximal
ideal of ${\cal{O}}_r$ and $X$ the germ defined by its maximal
minors. The normal module is given by
\begin{equation*}
{\mathcal{N}}_X\cong \dfrac{Mat_{(n+1,n)}({\cal{O}}_r)}{Im(g)}
\end{equation*}
where $g$ is the map defined as above.
\end{proposition}
\vspace{.1cm}
\begin{proposition}(\cite{FK}, pg.3997)
With the same notations of the previous lemma
\begin{equation*}
T_{X}^1=\dfrac{Mat_{(n+1,n)}({\cal{O}}_r)}{J(M)+Im(g)}.
\end{equation*}
\end{proposition}
\vspace{.1cm}

Our next goal is to write the theorem of finite determinacy for
$n\times p$ families of matrices.
In \cite{Bruce} and \cite{Farid}, the geometric characterization of
finitely determined germs of square matrices follows the ideas
given by Gaffney in \cite{wall}. In \cite{FK}, Anne Fr\"uhbis
characterizes finitely determined $(n+1)\times n$ matrices using
Artin Approximation Theorem. The proofs of the corresponding theorems (\ref{3.1}) and
(\ref{3.2}) are similar to the proofs the corresponding theorems for
Mather's groups (see \cite{FK} and \cite{tese}).

\begin{theorem}\label{3.1}[\textbf{Infinitesimal Criterion for Finite Determinacy}]
Let \linebreak $M\in Mat_{(n,p)}({\cal{O}}_r)$ and $k$ be a positive integer such that
\begin{equation*}
\mathcal{M}^{k+1}Mat_{(n,p)}({\cal{O}}_r)\subseteq
\mathcal{M}^2J(M)+{\mathcal{M}}Im(g).
\end{equation*}

Then, $M$ is $k$-finitely determined.
\end{theorem}

We observe that the action of $GL_{p}({\matC}) \times
GL_{n}({\matC})$ on $Mat_{(n,p)}(\matC)$ has $n$ orbits, given by $\Delta_i\backslash \Delta_{i-1}$, $1\leq i\leq n$, if $n\leq p$. These orbits determine a stratification of $Mat_{(n,p)}(\matC)$, which is a Whitney stratification (see \cite{EG}).

\begin{theorem}\label{3.2}\label{detfin}[\textbf{Geometric Criterion of Finite Determinacy}]
An element $M:\matC^r,0\longrightarrow Mat_{(n,p)}(\matC)$ is
$\cal{G}$-finitely determined if and only if $M$ is transverse to
the strata of the stratification of $Mat_{(n,p)}(\matC)$ outside the
origin.
\end{theorem}

\begin{corollary}
Let $M$ be a $n\times p$ matrix with entries in the maximal ideal of
${\cal{O}}_r$, defining an isolated singularity. Then $M$ is
$\cal{G}$-finitely determined.
\end{corollary}

\noindent \textbf{Proof.}
Since $M$ defines an isolated singularity, then $M(x)$ does not intercept $\Delta_i\backslash \Delta_{i-1}$
if $\Delta_i\backslash \Delta_{i-1}\neq \{0\}$. Then $M$ is
transverse to the strata of the stratification outside the origin, and by
the geometric criterium, $M$ is ${\cal{G}}$-finitely determined.

$\hfill \blacksquare$

\vspace{.3mm}

A different proof of this result can be found in \cite{FK}, pg.
3998. The result of Theorem \ref{3.1} can be generalized to
matrices on $Mat_{n,p}(\matR)$ with entries ${\cal{C}}^{\infty}$.

\section{The main result}

The purpose of this section is to prove the following theorem which
gives simple geometric and algebraic conditions characterizing
$\cal{G}$-finite determinacy of $n\times (n+1)$ matrices defining
isolated singularities $X=f^{-1}(0)$ with $f=(f_1,...,f_{n+1})$. As
before $f_i$ denotes the $n\times n$-minor of $M$ obtained by
removing the $i$-th-column of $M$, $i=1,\cdots,n+1$. 

We fix some notations:
\begin{itemize}
  \item[{\rm{a}}] $E_{ij}$ is the $n\times p$ matrix with $1$ at the
  $(i,j)$ position and zero otherwise;
  \item[b)] If $q\neq s$ and $\gamma\neq\nu$, we denote by 
    $$\Delta^{(q,s)}=\dfrac{\partial f_q}{\partial x_{\gamma}}\dfrac{\partial f_s}{\partial x_{\nu}}-\dfrac{\partial f_q}{\partial x_{\nu}}
  \dfrac{\partial f_s}{\partial x_{\gamma}}$$
   a $2\times 2$-minor of the jacobian matrix of $f$, where $1\leq r,\,s\leq n+1$ and $1\leq \nu,\,\gamma\leq r$.
  \item [d)] Let $J_f$ be the ideal generated by the $2\times
  2$ minors of the jacobian matrix of $f$, i. e., $J_f=\langle\Delta^{(q,s)}:\,1\leq\,q,\,s\leq n+1\rangle$ and let $$I_{\cal{G}}(M)=J_f+<f_1,...,f_{n+1}>.$$
  \item [e)] Let $M^j$ be the $n\times n$ matrix obtained removing the $j$-th column
  of $M$. Indicate by $\textrm{cof}^j(m_{su})$ the cofactor of the element
  $m_{su}$ in $M^j$.
  \item[f) ] Analogously, $M^{j}_{ki}$ is the $(n-1)\times
  (n-1)$ matrix  removing from $M^j$ the $k$-th row and the $i$-
  th column. We denote by $\textrm{cof}^j_{ki}(m_{su})$ the cofactor of the element
  $m_{su}$ in $M^j_{ki}$. When $n=1$, we consider that $\textrm{cof}^j_{ki}(m_{su})$=1.
\end{itemize}

\begin{theorem}\label{equi}
Let $M$ be the germ of an $n\times (n+1)$ matrix with entries in the
maximal ideal of $\,{\cal{O}}_r$. Then, the
following statements are equivalent:
\begin{itemize}
   \item [{\rm(a)}] M is ${\cal{G}}$-finitely determined and $X$  has isolated singularity;
   \item [{\rm(b)}]$X\cap V(J_f)=\{0\}$;
   \item [{\rm(c)}] $I_{\cal{G}}(M)\supseteq {\mathcal{M}}^{k}$, for some positive integer $k$.
\end{itemize}
\end{theorem}

\vspace{.2cm}

\noindent \textbf{Proof.}

\vspace{.2cm}

To prove that $(a)\Longrightarrow (b)$, let $\Delta_n\subset
Mat_{(n,n+1)}(\matC)$ be the set of singular matrices. We note that
$\Delta_n$ is an irreducible analytic variety. Also, if we consider
the diagram

$$\matC^r\stackrel{M}\longrightarrow Mat_{(n,n+1)}(\matC)\stackrel{\delta_n}\longrightarrow\matC^{n+1},$$
then $\Delta_n=\delta_n^{-1}(0)$. Now $M$ is ${\cal{G}}$-finitely determined and has isolated singularity, hence $f=\delta_n\circ M$ is a submersion away from zero. Therefore $X\cap V(J_f)=\{0\}$.

It follows from Hilbert-Nullstellensatz Theorem that
$(b)\Leftrightarrow (c)$. 


The proof $(c)\Longrightarrow(a)$ is harder. It is based on
Propositions (\ref{prop1}) and (\ref{gen}) in which we show that
the matrices $f_jE_{kl}$ and $J_fE_{kl}$ are on the
$\cal{G}$-tangent space to M, for all $1\leq l,\,j\leq n+1$ and
$1\leq k\leq n$. Then it follows from these conditions that 
$I_{\cal{G}}(M) Mat_{(n, n+1)}(\matC)\subseteq T_{\cal{G}}M$. Hence,
if $(c)$ holds, then $T_{\cal{G}}M\supset {\mathcal{M}}^{k}Mat_{(n,
(n+1))}$ and the result follows from Theorem (\ref{3.1}).

 $\hfill \blacksquare$ \vspace{.3cm}

\begin{proposition}\label{prop1}
Let $M$ be a $n\times (n+1)$ matrix with entries in the maximal
ideal of ${\cal{O}}_r$. Then $f_j E_{kl}\in T{\cal{G}} M$, for
$1\leq l,\,j\leq n+1$ and $1\leq k\leq n$.
\end{proposition}

\noindent \textbf{Proof.}

If $M=(m_{us})$, we can write each $f_j$, $j=1,...,n+1$, expanding
by the $k$-th-row of the matrix $M$ for any choice of $k$, $1\leq
k\leq n$, in the following way:
$$f_j=\sum_{{s\neq j}}m_{ks}\textrm{cof}^{j}(m_{ks}).$$

As before, we denote $C_{ls}(M)$, or $C_{ls}$ when the context is
clear, the matrix that has $l$-th column equal to $s$-th column of
$M$ and zeros in any other position, $1\leq l,s\leq n+1$.

For each $l$, we can consider the $n\times(n+1)$ matrix
$A_l=(a_{us})$ defined by
$$A_l=\sum_{{s\neq j}}\textrm{cof}^{j}(m_{ks})C_{ls}\in T{\cal{G}}M.$$
Then,
\begin{itemize}
   \item [i) ] $a_{us}=0$, if $s\neq l$;
   \item[ii) ] ${\displaystyle a_{kl}=\sum_{{s\neq
   l}}m_{ks}\textrm{cof}^{k}(m_{ks})= f_j}$;
   \item[iii) ] ${\displaystyle a_{ul}=\sum_{{s\neq
   j}}m_{us}\textrm{cof}^{j}(m_{ks})}$=0, for $u\neq k$, since this is
the determinant of a matrix that has two identical rows.
\end{itemize}
Therefore, $A_l=f_jE_{kl}$ and we get the result.

$\hfill \blacksquare$ \vspace{.3mm}

\begin{example}
Let $M=(m_{ij}(x))$ be a $3\times 4$ matrix with entries in the maximal ideal of ${\cal{O}}_r$, and we choose $l=2$ and $k=1$ in the previous lemma.  To verify that $f_1E_{12}\in\,T{\cal{G}}M,$  we write
$$A_2=\sum_{{s\neq 1}}\textrm{cof}^{1}(m_{1s})C_{2s}=\textrm{cof}^{1}(m_{12})C_{22}+\textrm{cof}^{1}(m_{13})C_{23}+\textrm{cof}^{1}(m_{14})C_{24}.$$
We note that 
\begin{itemize}
\item[i)] $cof^1(m_{12})=m_{23}m_{34}-m_{33}m_{24}$,
\item[ii)] $cof^1(m_{13})=-(m_{22}m_{34}-m_{32}m_{24})$,
\item[iii)] $cof^1(m_{14})=m_{22}m_{33}-m_{32}m_{23}$. 
\end{itemize}
Then
\begin{itemize}
\item[a)] $a_{12}=m_{12}\textrm{cof}^1(m_{12})+m_{13}\textrm{cof}^1(m_{13})+m_{14}\textrm{cof}^1(m_{14})=f_1$
\item[b)] $a_{22}=m_{22}\textrm{cof}^1(m_{12})+m_{23}\textrm{cof}^1(m_{13})+m_{24}\textrm{cof}^1(m_{14})=0$
\item[c)] $a_{32}=m_{32}\textrm{cof}^1(m_{12})+m_{33}\textrm{cof}^1(m_{13})+m_{34}\textrm{cof}^1(m_{14})=0$
\end{itemize}
Moreover if $t\neq 2$, $a_{it}=0$, $1\leq i\leq 3$. Therefore $f_1E_{12}$ belongs to the tangent space of $M$.
\end{example} 

\vspace{.3mm}

Our next gol is to proof that the $2\times 2$-minors of the jacobian matrix of $f$ are on the $\mathcal{G}$-tangent space of $M$. For this we need some preliminary lemmas.


\begin{lemma}\label{l4.1}
Let $M=(m_{us})$ be a $n\times (n+1)$ matrix with entries in the
maximal ideal of ${\cal{O}}_r$. Then,
$$\textrm{cof}^{\,j}(m_{il})=(-1)^{\alpha}\textrm{cof}^{\,\,l}(m_{ij}),$$
where $$\alpha=\left\{\begin{array}{ll}
               l-j+1,\mbox{ if } l<j \\
               l-j-1,\mbox{ if } l>j
             \end{array}\right.$$
$1\leq l,j\leq n+1$ and $1\leq i\leq n$.
\end{lemma}
\noindent \textbf{Proof.}

The proof follows directly from the expressions of $\textrm{cof}^j(m_{il})$
and $\textrm{cof}^l(m_{ij})$.

$\hfill \blacksquare$ \vspace{.5mm}

\begin{lemma}\label{l1}
We fix $j,\,l,$ and $\gamma$, such that $j\neq l$, $1\leq
j,\,l\leq n+1$, and $1\leq\gamma\leq r$. Then,
\begin{equation}\label{Ml1}
G_{jl}^{\gamma}=\dfrac{\partial f_j}{\partial x_{\gamma}}E_{kl}+(-1)^{l-j+1}\dfrac{\partial f_l}{\partial x_{\gamma}}E_{kj}\in T{\cal{G}},
\end{equation}
for   $1\leq k\leq n$.
\end{lemma}
\noindent \textbf{Proof.}

Without loss of generality, we let $l<j$. We derive the proof in four steps:

\noindent\textbf{\underline{Step 1:}} Let $A=(a_{qv})$ be the matrix
defined by
$$A=\dfrac{\partial M}{\partial x_{\gamma}}\textrm{cof}^j(m_{kl})+\sum_{i\neq j}\dfrac{\partial \textrm{cof}^j(m_{ki})}{\partial x_{\gamma}}C_{li}+\sum_{i\neq j,\, l}
(-1)^{\alpha}\dfrac{\partial m_{ki}}{\partial
x_{\gamma}}\left(\sum_{u\neq
k}R_{ku}\textrm{cof}^j_{ki}(m_{ul})\right),$$ where  $R_{us}$ is the
matrix that has $u$-th row equal to $s$-th row of $M$ with zeros in
any other position and
$$\alpha=\left\{\begin{array}{ll}
                                                             k+i,\mbox{ if } i<j \\
                                                            k+i-1, \mbox{ if }
                                                             i>j.
                                                           \end{array}\right.
$$
We now look at each entry of the matrix $A$. We deal in
$(i)$-$(iii)$ with the case $q=k$, and in $(iv)$ and $(v)$ with the
case $q\neq k$.
  \begin{itemize}
 \item[i) ] At the entry $(k,j)$, we have
  \begin{align*}
  a_{kj}=
  \dfrac{\partial m_{kj}}{\partial x_{\gamma}}\textrm{cof}^j(m_{kl})+\sum_{i\neq j}(-1)^{\alpha}
  \dfrac{\partial m_{ki}}{\partial x_{\gamma}}\left(\sum_{u\neq
  k}m_{uj}\textrm{cof}^j_{ki}(m_{ul})\right)
  \end{align*}
  Using lemma \ref{l4.1}, we can write
  \begin{align*}
  a_{kj}&=(-1)^{l-j+1}\left(\dfrac{\partial m_{kj}}{\partial x_{\gamma}}\textrm{cof}^l(m_{kj})+\sum_{i\neq j}(-1)^{\alpha}
  \dfrac{\partial m_{ki}}{\partial x_{\gamma}}\left(\sum_{u\neq
  k}m_{uj}\textrm{cof}^l_{ki}(m_{uj})\right)\right)\\&=(-1)^{l-j+1}\left(\dfrac{\partial f_l}{\partial x_{\gamma}}
  -\sum_{i\neq l}\dfrac{\partial \textrm{cof}^l(m_{ki})}{\partial
  x_{\gamma}}m_{ki}\right).
  \end{align*}
 \item[ii) ] At the entry $(k,l)$, we have
    \begin{align*}a_{kl}&=
  \dfrac{\partial m_{kl}}{\partial x_{\gamma}}\textrm{cof}^j(m_{kl})+\sum_{i\neq j}\dfrac{\partial \textrm{cof}^j(m_{ki})}{\partial x_{\gamma}}m_{ki}
  +\sum_{i\neq j,\,l}(-1)^{\alpha}\dfrac{\partial m_{ki}}{\partial x_{\gamma}}\left(\sum_{u\neq k}m_{ul}\textrm{cof}^j_{ki}(m_{ul})\right)\\
  &=\sum_{i\neq j}\dfrac{\partial \textrm{cof}^j(m_{ki})}{\partial x_{\gamma}}m_{ki}
  +\sum_{i\neq j}\dfrac{\partial m_{ki}}{\partial x_{\gamma}}\textrm{cof}^j(m_{ki})=\dfrac{\partial f_j}{\partial
  x_{\gamma}}.
  \end{align*}

\item[iii) ] If $t\neq l,\,j$, at the entry $(k,t)$ we have:
  \begin{align*}
  a_{kt}=\dfrac{\partial m_{kt}}{\partial x_{\gamma}}\textrm{cof}^j(m_{kl})+\sum_{i\neq j}(-1)^{\alpha}
  \dfrac{\partial m_{ki}}{\partial x_{\gamma}}\left(\sum_{u\neq k}m_{ut}\textrm{cof}^j_{ki}(m_{ul})\right).
  \end{align*}
Notice that for $i\neq t$,
$$\sum_{u\neq k}m_{ut}\textrm{cof}^j_{ki}(m_{ul})=0,$$
since it is the determinant of a $(n-1)\times (n-1)$ matrix with two
columns equal to the $l$-th column of the matrix $M$. Then,
$$a_{kt}=\dfrac{\partial m_{kt}}{\partial x_{\gamma}}\textrm{cof}^j(m_{kl})+(-1)^{k+t}\dfrac{\partial m_{kt}}{\partial x_{\gamma}}\left(\sum_{u\neq
  k}m_{ut}\textrm{cof}^j_{kt}(m_{ul})\right).$$

Now, it is not hard to verify that $a_{kt}=0$, for all $t\neq
l,\,j$.

\item[iv) ] Let $q\neq k$, then at the entry $(q,l)$ we have:
$$a_{ql}=\dfrac{\partial m_{ql}}{\partial x_{\gamma}}\textrm{cof}^j(m_{kl})+\sum_{i\neq j}
  \dfrac{\partial \textrm{cof}^j(m_{ki})}{\partial   x_{\gamma}}m_{qi}.$$
We can write
$\displaystyle{\textrm{cof}^j(m_{ki})=(-1)^{\beta}\sum_{t\neq i,\,j}
m_{qt}\textrm{cof}^j_{ki}(m_{qt})}$,
where
$$\beta=\left\{\begin{array}{ll}
                 k+i,\mbox{ if } i<j \\
                 k+i-1, \mbox{ if } i>j.
                     \end{array}\right.$$

Then,
\begin{equation*}\label{somadeltacof}
a_{ql}=\dfrac{\partial m_{ql}}{\partial
x_{\gamma}}\textrm{cof}^j(m_{kl})+\underbrace{ \sum_{i\neq
j}(-1)^{\beta}\sum_{t\neq i,\,j} \dfrac{\partial m_{qt}}{\partial
x_{\gamma}}\textrm{cof}^j_{ki}(m_{qt})m_{qi}}_{S_1}+\underbrace{\sum_{i\neq
j}(-1)^{\beta}m_{qt}\dfrac{\partial
\textrm{cof}^j_{ki}(m_{qt})}{\partial x_{\gamma}} m_{qi}}_{S_2}.
\end{equation*}

With similar arguments as in the previous steps, we can show that
$$S_1=\sum_{t\neq j} -\dfrac{\partial m_{qt}}{\partial
x_{\gamma}}\textrm{cof}^j(m_{kt}),$$ and $S_2=0$.

Therefore,
$$a_{ql}= \sum_{t\neq j,\,l} -\dfrac{\partial m_{qt}}{\partial
x_{\gamma}}\textrm{cof}^j(m_{kt})$$

\item[v) ] The entry $(q,v)$, $v\neq l$ and $q\neq k$,
$$a_{qv}=\dfrac{\partial m_{qv}}{\partial x_{\gamma}}\textrm{cof}^j(m_{kl}).$$

\end{itemize}

\noindent\textbf{\underline{Step 2:}} We consider the $n\times
(n+1)$ matrix, $B=(b_{us})$, given by:
$$B=(-1)^{l-j+1}\sum_{i\neq l}\dfrac{\partial \textrm{cof}^l(m_{ki})}{\partial x_{\gamma}}C_{ji}.$$

Then,
\begin{itemize}
  \item[i) ] $b_{us}=0$, if $s\neq j$;
  \item[ ii)] $b_{uj}=\displaystyle{(-1)^{l-j+1}\sum_{i\neq l}\dfrac{\partial \textrm{cof}^l(m_{ki})}{\partial x_r}m_{ui}}$
\end{itemize}

\noindent\textbf{\underline{Step 3:}} We consider the $n\times
(n+1)$ matrix, $C=(c_{qv})$ given by $C=A+B$. Then,
  \begin{itemize}
    \item[i) ]$\displaystyle{c_{kj}=(-1)^{l-j+1}\dfrac{\partial f_l}{\partial x_{\gamma}}}$;
    \item[ii) ]$\displaystyle{c_{kl}=\dfrac{\partial f_j}{\partial x_{\gamma}}};$
    \item[iii) ]$c_{kt}=0$, for $t\neq l,j$;
    \item[iv) ]$\displaystyle{c_{ql}= -\sum_{t\neq j,\,l} \dfrac{\partial m_{qt}}{\partial
    x_{\gamma}}\textrm{cof}^j(m_{kt})}$, for $q\neq k$.
    \item[v) ] $\displaystyle{c_{qj}=\dfrac{\partial m_{qj}}{\partial x_{\gamma}}\textrm{cof}^j(m_{kl})+(-1)^{l-j+1}
    \sum_{i\neq l}\dfrac{\partial \textrm{cof}^j(m_{ki})}{\partial
    x_r}m_{ui}}$, for $q\neq k$.

    As in item \textit{
    (iv)} of the step $1$, it is possible to show that
$$c_{qj}= (-1)^{l-j+1}\sum_{t\neq j,\,l} \dfrac{\partial m_{qt}}{\partial
    x_{\gamma}}\textrm{cof}^{\,l}(m_{kj}).$$
    \item[vi) ] $c_{qv}=\dfrac{\partial m_{qv}}{\partial x_{\gamma}}\textrm{cof}^j(m_{kl})$, $\forall\,\,v\neq l,j$ and $q\neq
    k$.
  \end{itemize}

\noindent\textbf{\underline{Step 4:}} For each $q\neq k$, $1\leq
q\leq n$, we consider the matrices
$$D_q=\sum_{i\neq j,\,l}\dfrac{\partial m_{qi}}{\partial x_{\gamma}}\left(\sum_{u\neq k}(-1)^{\mu}R_{qu}\textrm{cof}^j_{ui}(m_{kl})\right),$$
where
$$\mu=\left\{\begin{array}{l}
                   u+i-1,\mbox{ if } j<i \mbox{ and } u>k \mbox{ or } j>i \mbox{ and } u>k\\
                   u+i,\mbox{ if } j<i \mbox{ and } u<k  \mbox{ or } j>i \mbox{ and }
                   u<k.
                   \end{array}\right.$$

We look at each entry of the matrix $D$:
\begin{itemize}
\item [i) ] $d_{ql}=\displaystyle{\sum_{i\neq j,\,l}\dfrac{\partial m_{qi}}{\partial x_{\gamma}}\left(\sum_{u\neq k}(-1)^{\gamma}m_{ul}
\textrm{cof}^j_{ui}(m_{kl})\right)}.$

Using lemma \ref{l4.1} and the expression of the
$\textrm{cof}^j(m_{ki})$ it is not difficult to see that
$$\sum_{u\neq k}(-1)^{\gamma}m_{ul}\textrm{cof}^j_{ui}(m_{kl})=\textrm{cof}^j(m_{ki}).$$

Thus,
$$d_{ql}=\sum_{i\neq j,\,l}\dfrac{\partial m_{qi}}{\partial x_{\gamma}}\textrm{cof}^j(m_{ki}).$$

\item [ii) ] $d_{qj}=\displaystyle{\sum_{i\neq j,\,l}\dfrac{\partial m_{qi}}{\partial x_{\gamma}}\left(\sum_{u\neq k}(-1)^{\gamma}m_{uj}
\textrm{cof}^j_{ui}(m_{kl})\right)}.$ Again, using lemma \ref{l4.1}
and the expression of the $\textrm{cof}^l(m_{kj})$ we have
$$d_{qj}=\displaystyle{(-1)^{j-l-1}\sum_{i\neq j,\,l}\dfrac{\partial m_{qi}}{\partial x_{\gamma}}
\textrm{cof}^{\,l}(m_{kj})}$$

\item[iii) ] For $t\neq j,\,l$, we have
$$d_{qt}=\displaystyle{\sum_{i\neq j,\,l}\dfrac{\partial
m_{qi}}{\partial x_{\gamma}}\left(\sum_{u\neq k}(-1)^{\gamma}m_{ut}
\textrm{cof}^j_{ui}(m_{kl})\right)}.$$

Note that for $i\neq t$, we have
$$\sum_{u\neq k}(-1)^{\gamma}m_{ut}
\textrm{cof}^j_{ui}(m_{kl})=0$$ since this is the expression of the
determinant of a matrix that has two equal columns.

Then,
$$d_{qt}=\displaystyle{\dfrac{\partial
m_{qt}}{\partial x_{\gamma}}\left(\sum_{u\neq k}(-1)^{\gamma}m_{ut}
\textrm{cof}^j_{ut}(m_{kl})\right)}.$$

Using lemma \ref{l4.1} and the expression of the
$\textrm{cof}^j(m_{kl})$ we have
$$d_{qt}=-\displaystyle{\dfrac{\partial
m_{qt}}{\partial x_{\gamma}}\textrm{cof}^j(m_{kl})}$$

\item[i) ] $d_{sv}=0$, for $s\neq q$;
  \end{itemize}

To conclude the proof, it suffices to consider the matrix
$E=(e_{qv})$ given by $\displaystyle{E=C+\sum_{q\neq k}D_q}$.

$\hfill \blacksquare$ \vspace{.5cm}

\begin{example}
\label{ex1}
Let $M=(m_{ij})$ be a $2\times 3$ matrix with entries in the maximal ideal of ${\cal{O}}_4$. Then using the previous result, we will verify
$$\dfrac{\partial f_2}{\partial x_{\gamma}}E_{11}+\dfrac{\partial f_1}{\partial x_{\gamma}}E_{12}$$
belongs to $T{\cal{G}}M$. In fact, we first consider the matrix 
$$A=\dfrac{\partial M}{\partial x_{\gamma}}\textrm{cof}^2(m_{11})+\dfrac{\partial \textrm{cof}^2(m_{11})}{\partial x_{\gamma}}C_{11}+\dfrac{\partial \textrm{cof}^2(m_{13})}{\partial x_{\gamma}}C_{13}-\dfrac{\partial m_{13}}{\partial x_{\gamma}}R_{12}$$
as in the step 1 of the Lemma \ref{l1}. Then, we look at each entry of $A$:
\begin{itemize}

\item[a)] $a_{11}=m_{23}\dfrac{\partial m_{11}}{\partial x_{\gamma}}+m_{11}\dfrac{\partial m_{23}}{\partial x_{\gamma}}-m_{13}\dfrac{\partial m_{21}}{\partial x_{\gamma}}-m_{21}\dfrac{\partial m_{13}}{\partial x_{\gamma}}=\dfrac{\partial f_2}{\partial x_{\gamma}}$
\item[b)]$a_{12}=m_{23}\dfrac{\partial m_{12}}{\partial x_{\gamma}}-m_{22}\dfrac{\partial m_{13}}{\partial x_{\gamma}}=\dfrac{\partial f_1}{\partial x_{\gamma}}-\left(m_{12}\dfrac{\partial m_{23}}{\partial x_{\gamma}}+m_{13}\dfrac{\partial m_{22}}{\partial x_{\gamma}}\right)$
\item[c)] $a_{13}=m_{23}\dfrac{\partial m_{13}}{\partial x_{\gamma}}-m_{23}\dfrac{\partial m_{13}}{\partial x_{\gamma}}=0$
\item[d)] $a_{21}=m_{23}\dfrac{\partial m_{21}}{\partial x_{\gamma}}+m_{21}\dfrac{\partial m_{23}}{\partial x_{\gamma}}-m_{23}\dfrac{\partial m_{21}}{\partial x_{\gamma}}=m_{21}\dfrac{\partial m_{23}}{\partial x_{\gamma}}$
\item[e)] $a_{22}=m_{23}\dfrac{\partial m_{22}}{\partial x_{\gamma}}$
\item[f)] $a_{23}=m_{23}\dfrac{\partial m_{23}}{\partial x_{\gamma}}$
\end{itemize}

Now, we define the matrix
$$B=\dfrac{\partial \textrm{cof}^1(m_{12})}{\partial x_{\gamma}}C_{22}+\dfrac{\partial \textrm{cof}^1(m_{13})}{\partial x_{\gamma}}C_{23}=\left(\begin{array}{ccc}
0 & m_{12}\dfrac{\partial m_{23}}{\partial x_{\gamma}}-m_{13}\dfrac{\partial m_{22}}{\partial x_{\gamma}}& 0 \\
0 & m_{22}\dfrac{\partial m_{23}}{\partial x_{\gamma}}-m_{23}\dfrac{\partial m_{22}}{\partial x_{\gamma}}& 0 \\
                           \end{array}
                         \right),$$
according to the step 2 of the Lemma \ref{l1}.
Then 

$$A+B=\left(\begin{array}{ccc}
\dfrac{\partial f_2}{x_{\gamma}} & \dfrac{\partial f_1}{x_{\gamma}}& 0 \\
m_{21}\dfrac{\partial m_{23}}{x_{\gamma}}& m_{22}\dfrac{\partial m_{23}}{x_{\gamma}}& m_{23}\dfrac{\partial m_{23}}{x_{\gamma}} \\
                           \end{array}
                         \right).$$
Finally, by the step 4, we have the matrix
$$D_{2}=-\dfrac{\partial m_{23}}{\partial x_{\gamma}}R_{22}=-\left(
                           \begin{array}{ccc}
                             0 & 0 & 0 \\
                             m_{21}\dfrac{\partial m_{23}}{\partial x_{\gamma}} & m_{22}\dfrac{\partial m_{23}}{\partial x_{\gamma}} & m_{23}\dfrac{\partial m_{23}}{\partial x_{\gamma}} \\
                           \end{array}
                         \right),$$
and the result follows from adding $A+B$ with $D_2$.
\end{example}

\begin{proposition}\label{gen}
The matrices $\Delta_{(j,\,t)}E_{kl}$ belong to the
${\cal{G}}$-tangent space of the matrix $M$, for $1\leq k\leq n$ and
$1\leq j,\,t,\,l\leq n+1$, $j\neq t$.
\end{proposition}
\noindent \textbf{Proof.}

We show that $\Delta^{(j,\,t)}E_{kl}$ are obtained
using the matrices of the previous lemma. Let us consider two cases:
\vspace{.5mm}

\noindent\textbf{i)}  If $j\neq l$ and $t\neq l$, then
\begin{align*}
&\Delta^{(j,\,t)}E_{kl}=(-1)^{l-j}\dfrac{\partial f_l}{\partial x_{\nu}}\left( \dfrac{\partial f_t}{\partial x_{\gamma}}E_{kj}+(-1)^{j-t+1}\dfrac{\partial f_j}{\partial x_{\gamma}}E_{kt}\right)+\\ &+\dfrac{\partial f_t}{\partial x_{\gamma}}\left(\dfrac{\partial f_j}{\partial x_{\nu}}E_{kl}+(-1)^{l-j+1}\dfrac{\partial f_l}{\partial x_{\nu}}E_{kj}\right)-\dfrac{\partial f_j}{\partial x_{\gamma}}\left(\dfrac{\partial f_{t}}{\partial x_{\nu}}E_{kl}+(-1)^{l-t+1}\dfrac{\partial f_l}{\partial x_{\nu}}E_{kt}\right),
\end{align*}
then $\Delta^{(j,\,t)}E_{kl}\in T{\cal{G}}M.$
\vspace{.5mm}

\noindent\textbf{ii)} If $j= l$, then
\begin{align*}
&\Delta^{(l,\,t)}E_{kl}=\\
&=\dfrac{\partial f_l}{\partial x_{\gamma}}\left(\dfrac{\partial f_t}{\partial x_{\nu}}E_{kl}+(-1)^{l-t+1}\dfrac{\partial f_l}{\partial x_{\nu}}E_{kt}\right)-\dfrac{\partial f_l}{\partial x_{\nu}}\left(\dfrac{\partial f_t}{\partial x_{\gamma}}E_{kl}+(-1)^{l-t+1}\dfrac{\partial f_l}{\partial x_{\gamma}}E_{kt}\right)
\end{align*}
    Therefore, $\Delta^{(l,\,t)}E_{kl}\in T{\cal{G}}M.$

$\hfill \blacksquare$ \vspace{.5cm}


%

\section{${\cal{G}}$-Topological Equivalence of Matrices}

As an application of the results of the previous section we study the
${\cal{G}}$-topological triviality of families of matrices.

We concentrate our study on $n\times (n+1)$ matrices $M$ with
entries in ${\cal{M}}_r$. We denote by $f=(f_1,\ldots,f_{n+1})$ the
ideal generated by its  maximal minors, $X$ the  variety defined by
$f$ and by ${\cal{O}}_r^0$ the ring of germs at the origin of
continuous functions of $\matK^r\longrightarrow\matK$, where
$\matK=\matR$ or $\matC$. The results are applied to
deformations of germs of weighted homogeneous matrices.

\begin{definition}
Two  matrices $M,\,N \in Mat_{(n,n+1)}({\cal{O}}_r)$ are topologically equivalent (or
$C^0$-${\cal{G}}$-equivalent) if there exist a germ of homeomorphism
$\phi:(\matK^r,0)\longrightarrow (\matK^r,0)$ and invertible
matrices $A\in GL_{n}({\cal{O}}_r^0)$ and $B\in
GL_{n+1}({\cal{O}}_r^0)$ such that $M=A^{-1}(N\circ \phi)B$.
\end{definition}

A control function $\rho:\matK^r\longrightarrow \matR$ is a non
negative function that satisfies the following condition:
\begin{itemize}
  \item[i)] $\rho(0)=0$ and there exist constants $c>0$ and
  $\alpha>0$ such that $\rho(x)\geq c|x|^{\alpha}$ (i. e., $\rho$ satisfies a Lojasiewicz condition).
\end{itemize}

\begin{definition}
A matrix $M=(m_{ij})\in Mat_{(n,n+1)}({\cal{O}}_r)$ is
$k$-$C^0$-${\cal{G}}$-determined, if for every matrix $N=(n_{ij})\in Mat_{(n,n+1)}({\cal{O}}_r)$ such that
$j^kM(0)=j^kN(0)$ \footnote{$j^kM(0)$ denotes the $n\times(n+1)$ matrix whose entries are the $k$-jets of $m_{ij}$ at zero}, $1\leq i\leq n$ and $1\leq j \leq n+1$, $N$ is $C^0$-${\cal{G}}$-equivalent to $M$.
\end{definition}

\begin{definition}\label{Trivial}
A one parameter deformation $M\in Mat_{(n,n+1)}({\cal{O}}_{r+1})$ of $M_0\in Mat_{(n,n+1)}({\cal{O}}_r)$ is $C^0-{\cal{G}}$- trivial (or ${\cal{G}}$-trivial) if there exists a homeomorphism
\begin{align*}
\Phi:(\matK^r\times\matK,0)&\longrightarrow(\matK^r\times\matK,0)\\
(x,t)&\longmapsto(\phi(x,t),t)
\end{align*}
such that $\Phi(x,0)=(x,0)$, $\phi(0,t)=0$ and families of matrices $A\in GL_{n}({\cal{O}}_{r+1}^0)$ and $B\in
GL_{n+1}({\cal{O}}_{r+1}^0)$ such that $A(x,0)=I_{n}$, $B(x,0)=I_{n+1}$
and $M_0=A^{-1}(M\circ\Phi)B$.

\end{definition}

\begin{proposition}\label{thom}
Let $M_0$ be a $n\times (n+1)$ matrix, defining a
codimension $2$ Cohen- Macaulay isolated singularity and $M$ a
deformation of $M_0$. Suppose that there exists a control function
$\rho$ such that
$$\rho^2\dfrac{\partial M}{\partial
t}=\sum_{i=1}^{r}\xi_{i}\dfrac{\partial M}{\partial x_{i}}+\sum_{l,k=1}^{n+1}L_{lk}C_{lk}(M)+\sum_{r,s=1}^{n}S_{rs}R_{rs}(M),$$
with $\xi_i(x,t),\, L_{lk}(x,t),S_{rs}(x,t)\in {\cal{O}}_{r+1}$, satisfying the conditions
\begin{align}\label{*}
&\dfrac{|\xi_i(x,t)|}{\rho^2(x,t)}\leq
C_1|x|,\,\,\,\dfrac{|L_{lk}(x,t)|}{\rho^2(x,t)}\leq C_2|x|,\,\,\,
,\dfrac{|S_{rs}(x,t)|}{\rho^2(x,t)}\leq C_3|x|.
\end{align}
with $\xi_i(0,t)=0$, $L(x,0)=Id_{n+1}$ e $S(x,0)=Id_{n}$.
Then the family $M$ is $C^0$-${\cal{G}}$-trivial.

\end{proposition}
\noindent \textbf{Proof.}

To get topological triviality of $M$ we construct continuous vector
fields $X$ in $\matK^r\times\matK$, $W$ in $\matK^{n^2}\times\matK$ and $Z$ in $\matK^{(n+1)^2}\times\matK$, lifting $\dfrac{\partial}{\partial t}$. The condition (\ref{*}) ensures the uniqueness of the corresponding flow, which gives the family of homeomorphism trivializing $M$.

We start the proof defining $\Phi$. Let $X(x,t)$ be the vector
field in $\matK^r\times\matK,0$ defined by
$$\left\{\begin{array}{ll}
    \dfrac{\partial}{\partial t}-\displaystyle{\sum_{i=1}^{r}}\dfrac{\xi_i(x,t)}{\rho^2(x,t)}\dfrac{\partial}{\partial x_i}, \mbox{ if } x\neq 0\\
    \dfrac{\partial}{\partial t}, \mbox{ if } x=0.
  \end{array}\right.
$$

For each $j$, $X_j$ denotes the $j$-th component of $X$. The vector
field $X$ is real analytic along
$(\matK^r\times\matK)\backslash(\{0\}\times \matK)$. Furthermore,
$$|X_j(x,t)|=\dfrac{|\xi_j(x,t)|}{\rho^2(x,t)}\leq C_1|x|,$$
$1\leq j\leq r$, so that the vector field $X(x,t)$ satisfies a Lipschitz
condition along $\{0\}\times\matK$. It follows from \cite{kuo} that
this vector field is locally integrable. For more details on the complex
case, see \cite{ruas}. We indicate $\Phi(x,t)$ the corresponding flow.

Now we want to find matrices $A$ and $B$ satisfying the conditions of Definition (\ref{Trivial}), that is, $A\in GL_{n}({\cal{O}}_{r+1}^0)$ and $B\in
GL_{n+1}({\cal{O}}_{r+1}^0)$ such that $A(x,0)=I_{n}$, $B(x,0)=I_{n+1}$
and $AM_0B^{-1}=M\circ\Phi$. 
 
By definition of $\Phi$ and the hypothesis, we have

\begin{align*}
\dfrac{\partial (m_{ij}\circ\Phi)}{\partial t}(x,t)&=-\sum_{i=1}^{r}\xi_{i}(x,t)\dfrac{\partial m_{ij}}{\partial x_{i}}\Phi(x,t)+\dfrac{\partial m_{ij}}{\partial t}(x,t)=\\
&=-\left(\sum_{k=1}^{n+1}\dfrac{L_{jk}(\Phi(x,t))}{\rho^2(\Phi(x,t))}m_{ik}(\Phi(x,t))+\sum_{s=1}^{n}\dfrac{S_{is}(\Phi(x,t))}{\rho^2(\Phi(x,t))}m_{sj}(\Phi(x,t))\right).
\end{align*}
where $1\leq i\leq n$ and $1\leq j \leq n+1$. 

Then we want to find matrices $A$ and $C=B^{-1}$ such that
\begin{equation}\label{sismat1}
\dfrac{\partial }{\partial
t}(AM_0C)_{ij}=-\left(\sum_{k=1}^{n+1}\dfrac{L_{jk}(\Phi(x,t))}{\rho^2(\Phi(x,t))}(AM_0C)_{ik}+\sum_{s=1}^{n}\dfrac{S_{is}(\Phi(x,t))}{\rho^2(\Phi(x,t))}(AM_0C)_{sj}\right)
\end{equation}

To solve the system (\ref{sismat1}) it is sufficient to solve the following two systems of differential equations
\begin{align}
&\dfrac{\partial a_{ij}}{\partial t}=-\sum_{k=1}^{n+1}\dfrac{L_{jk}(\Phi(x,t))}{\rho^2(\Phi(x,t))}a_{ik}\label{si1}\\
&\dfrac{\partial c_{ij}}{\partial t}=-\sum_{s=1}^{n}\dfrac{S_{is}(\Phi(x,t))}{\rho^2(\Phi(x,t))}c_{sj}\label{si2}
\end{align} 
with the same initial conditions. In fact, $M_0$ is independent of $t$ and the original system appears multiplying (\ref{si1}) to the right by $M_0C$,
(\ref{si2}) to the left by $AM_0$ and  adding the resulting systems.

Consider the following vector field $W=(w_{ij})$ on
$\matK^{(n+1)^2}\times\matK$
\begin{align}\label{4}
&w_{ij}(x,y,t)=\dfrac{\partial}{\partial t}+\sum_{k=1}^{n+1}
\dfrac{L_{jk}(\Phi(x,t))}{\rho^2(\Phi(x,t))}y_{ik}.
\end{align}
In (\ref{4}) we use $y_{ik}$ to denote the variables of
{$\matK^{(n+1)^2}\times\matK$} where $1\leq j,\,k\leq n$. By
hypothesis,
\begin{align*}
&\dfrac{L_{ij}(\Phi(x,t))}{\rho^2(\Phi(x,t))}\leq C_2|x|
\end{align*}
It follows from \cite{ruas} that the vector field is integrable.

With similar arguments, we obtain that the vector field $Z=(z_{ij})$ in
$\matK^{n^2}\times\matK$ defined by
\begin{align*}
&z_{ij}(x,y,t)=\dfrac{\partial}{\partial t}+\sum_{s=1}^{n}
\dfrac{S_{is}(\Phi(x,t))}{\rho^2(\Phi(x,t))}y_{sj}
\end{align*}
is integrable.

$\hfill \blacksquare$

\vspace{0.3cm}
\section{Deformations of Weighted Homogeneous\newline Families of Matrices}

Our next goal is to prove the topological triviality theorem for
deformations of  $\cal{G}$-finitely 
determined weighted homogeneous matrices 
$M\in Mat_{(n,n+1)}({\cal{O}}_r)$. We will show
that for deformations of degree greater than the maximum weighted degree
of the entries of $M$, the hypothesis of Proposition
(\ref{thom}) holds, hence they are topologically $\cal{G}$-trivial.

\begin{definition}
Given a set of weights $a=(a_1,...,a_r)\in \matN^r$, for any monomial
$x^{\alpha}=x_1^{{\alpha}_1}x_2^{{\alpha}_2}...x_r^{{\alpha}_r}$,
we define the $a$-filtration of $x^{\alpha}$ by
$$fil(x^{\alpha})=\sum^{r}_{i=1}a_i{\alpha}_i.$$
\end{definition}

We can define a filtration on the ring ${\cal{O}}_r$, as follows
$$fil(f)=\inf_{\alpha}\left\{fil(x^{\alpha})|\dfrac{\partial^{|\alpha|} f}{\partial x^{\alpha}}(0)\neq
0\right\},$$ for all germ $f\in{\cal{O}}_r$, where $|\alpha|=\alpha_1+\alpha_2+\ldots+\alpha_r$. 

We extend the filtration to the submodule $\Theta_X$ of germs of vector fields tangent to $X$, defining $a\left(\dfrac{\partial}{\partial x_j}\right)=-a_j$ for all $j=1,\ldots, n$ so given $\xi=\sum_{i=1}^n\xi_j\dfrac{\partial}{\partial x_j}\in\Theta_X$, then $fil(\xi)=\inf_j\{fil(\xi_j)-w_j\}$.

We can extend this
definition to the ring of $1$-parameter families of germs on $r$
variables, putting
$$fil(x^{\alpha}t^{\beta})=fil(x^{\alpha}).$$

Given a matrix $M\in Mat_{(n,n+1)}({\cal{O}}_r)$, we define
$fil(M)={\cal{D}}=(d_{ij})$, where $d_{ij}=fil(m_{ij})$ for $1\leq i\leq
n$, $1\leq j\leq n+1$.

\vspace{0.3cm}

\begin{definition}
A matrix $M\in Mat_{(n,n+1)}({\cal{O}}_r)$ is called
weighted homogeneous of type $({\cal{D}};a)\in
Mat_{(n,n+1)}(\matN)\times \matN^r$, if
\begin{itemize}
   \item[i)] $fil(m_{ij})=d_{ij}$ with respect to $\,\,a=(a_1,...,a_r)$;
   \item[ii)] The following relations are verified
   $$d_{ij}-d_{ik}=d_{lj}-d_{lk} \mbox{ for all } 1\leq i,\,l\leq n,\,\,1\leq j,\,k\leq n+1.$$
\end{itemize}
\end{definition}

Let $M$ be a $n\times (n+1)$ weighted homogeneous matrix of type
$({\cal{D}};a)$. Let $f=(f_1,\ldots,f_{n+1})$ be the ideal generated by the maximal minors of $M$. The index of $f_u$ indicates the column removed from $M$ to compute the minor. Then it is immediate that $f$ is weighted
homogeneous of type  $(D_1,\ldots,D_{n+1};a)$, where
$fil(f_u)=D_u$. The converse is
also true, that is, if $f\in{\cal{O}}_r$ is a codimension $2$
Cohen-Macaulay ideal generated by weighted homogeneous polynomials
with respect to some set of weights $a$, then there exists a weighted
homogeneous presentation matrix $M$ of $f$ of type $({\cal{D}},a)$ for some
${\cal{D}}\in Mat_{(n,n+1)}(\matN)$ (see \cite{FK} for a proof).

Let $k_1={\displaystyle {l.c.m}\{D_u|1\leq u\leq n+1\}}$ and
$\beta_u=k_1/D_u$. We define
$${\cal{N}}_{\cal{H}}M={\displaystyle\sum_{j=1}^{n+1}|f_j|^{2\beta_j}}.$$
Note that ${\cal{N}}_{\cal{H}}M$ is a weighted homogeneous control function of
type $(2k_1;a)$.

We note that each $\dfrac{\partial f_q}{\partial x_{\gamma}}$
is weighted homogeneous of type $(D_{q}-a_{\gamma};a)$ and for each minor $\Delta^{(q,s)}$
of the Jacobian matrix of $f$, there exists an integer $D_{ij}$ such that 
$\Delta^{(q,s)}$ is weighted homogeneous of type $(D_{qs};a)$. 
For $k_2=l.c.m.(D_{ij})$ and $\alpha_{qs}=k_2/D_{qs}$ we define
$${\cal{N}}_{\cal{R}}M={\displaystyle\sum_{(q,\,s)}
|\Delta^{(q,s)}|^{2\alpha_{qs}}}.$$ This is a weighted homogeneous
function of type $(2k_2;a)$.

Let $K=l.c.m.\{k_1,k_2\}$ and $c_i=K/k_i$. We define
$${\cal{N}_{\cal{G}}}M={\cal{N}}_{\cal{R}}^{c_1}M+{\cal{N}}_{\cal{H}}^{c_2}M.$$
Thus, ${\cal{N}_{\cal{G}}}M$ is weighted homogeneous of type
$(2K;a)$.

Let $M_t(x)=M_0(x)+t\theta(x)$ be a deformation of $M_0$ with
$t\in[0,1]$, and $d_{max}=\displaystyle{\max_{ij}}\{d_{ij}\}$ and let $F$ denote the maximal minors of $M_t$.

We define the control function
${\cal{N}}_{\cal{R}}M_t={\displaystyle\sum_{i\in
I}|\Delta^{(q,s)}_t|^{2\alpha_{ij}}}$, where $\Delta^{(q,s)}_{t}$ are the $2\times
2$ minors of the jacobian matrix of $f_t$ and the $\alpha_{ij}$ are the
same as above. If $fil(\Theta_{ij})\geq d_{max}+1$, for all $i,j$, then there exist 
constants $C_1$ and $C_2$ such that $C_1{\cal{N}}_{\cal{R}}M_0\leq{\cal{N}}_{\cal{R}}M_t\leq C_2{\cal{N}}_{\cal{R}}M_0$. If $\Theta=(\theta_{ij})$ and $fil(\theta_{ij})= d_{max}$, 
then this condition also holds to $t$ small enough (see \cite{saia}).

We can define similarly, the control
${\cal{N}}_{\cal{H}}M_t={\displaystyle\sum_{i=1}^{n+1}|F_{i}|^{2\beta_i}}$,
where each $\beta_{i}$ is obtained as above and $F_{i}$ are the
$i$th component of $F$.

We will prove the topological triviality theorem in two  parts.
First, we show the result for the group ${\cal{H}}$. Then we extend
the proof to the case of the group ${\cal{G}}$.

\begin{proposition}\label{trivH}
Let $M_0\in Mat_{(n,n+1)}({\cal{O}}_r)$ be
the germ of a weighted homogeneous matrix of type $(D,a)\in
Mat_{(n,n+1)}(\matN)\times \matN^r $,
satisfying the condition ${\cal{N}}_{\cal{H}}(M_0(x))\geq
c|x|^{\alpha}$ for constants $c$ and $\alpha$. Then,
\begin{itemize}
\item[ i) ] Deformations $M(x,t)=M_0(x)+t\Theta(x)$ of $M_0$ with
$fil(\theta_{ij})\geq d_{max}+1$, $t\in[0,1]$, $\forall\, i,j$,  and
$d_{max}=\displaystyle{\max_{ij}}\{d_{ij}\}$, are
${\cal{C}}^{0}-{\cal{H}}-$trivial.

\item[ ii) ] Moreover, for $t$ small enough, deformations $M(x,t)=M_0(x)+t\Theta(x)$ of $M_0$ with
$fil(\theta_{ij})= d_{\max}$, $\forall\, i,j$, are
${\cal{C}}^{0}-{\cal{H}}-$trivial.
\end{itemize}
\end{proposition}
\noindent \textbf{Proof.}

Let $M(x,t)=M_0(x)+t\Theta(x)$ be a deformation
of $M_0$ with $fil(\theta_{ij})\geq d_{max}+1$, $t\in[0,1]$, $\forall\, i,j$. To obtain the ${\cal{C}}^0-{\cal{H}}-$triviality we find ${\cal{C}}^0$-map germs $L:\matK^{(n+1)^2}\times\matK\to \matK$ 
such that
$$\dfrac{\partial M}{\partial t}=\sum_{p=1}^{n+1}\sum_{s=1}^{n+1}L_{ps}(x,t)C_{ps}(M),$$
where $C_{rl}=C_{rl}(M)$ is defined in Section 2.

Using the Proposition \ref{prop1} we have

\begin{equation*}
f_{tj}\dfrac{\partial M}{\partial t}=\sum_{k=1}^{n}\sum_{l=1}^{n+1}\theta_{kl}\left(\sum_{i\ne j}\textrm{cof}^j(m_{ki})C_{li}(M)\right)
\end{equation*}
for $1\leq j\leq n+1$. 

Then, multiplying this equation by $|f_{tj}|^{2(\beta_j-1)}\overline{f}_{tj}$ and adding in $j$, we get
\begin{align*}
{\cal{N}}_{\cal{H}}M\dfrac{\partial M}{\partial t}&=\sum_{j=1}^{n+1}\left(\sum_{k=1}^{n}\sum_{l=1}^{n+1}\theta_{kl}\sum_{i\neq j}|f_{tj}|^{2(\beta_j-1)}\overline{f}_{tj}\textrm{cof}^j(m_{ki})C_{li}(M)\right)=\\
&=\sum_{l=1}^{n+1}\sum_{i=1}^{n+1}\left(\sum_{k=1}^{n}\theta_{kl}\sum_{j\neq i}|f_{tj}|^{2(\beta_j-1)}\overline{f}_{tj}\textrm{cof}^j(m_{ki})\right)C_{li}(M),
\end{align*}

We define,
$$L_{li}(x,t)=\sum_{k=1}^{n}\theta_{kl}\sum_{i\neq j}|f_{tj}|^{2(\beta_j-1)}\overline{f}_{tj}\textrm{cof}^j(m_{ki}).$$
Then,
$$\dfrac{\partial M}{\partial t}=\sum_{l=1}^{n+1}\sum_{i=1}^{n+1}\dfrac{L_{li}(x,t)}{{\cal{N}}_{\cal{H}}M}C_{li}(M).$$

Now
\begin{itemize}
   \item [i)] $\displaystyle{fil\left(\frac{\partial m_{ij}}{\partial t}\right)=fil\left(\theta_{ij}\right)\geq d_{max}+1}$;
   \item [ii)] $fil(f_{tj}^{2(\beta_j-1)}\overline{f}_{tj})=2k_1-D_j$;
   \item[iii)] $\displaystyle{fil(L_{ji})\geq 2k_1-D_{j}+d_{max}+1+D_{j}-d_{ki}\geq 2k_1+1}$, for all $k$.
\end{itemize}

Then, for each $1\leq j,\,i\leq n+1$,
$\dfrac{L_{ji}(x,t)}{{\cal{N}}_{\cal{H}}M_t}\leq C|x|$ and,
therefore, as in the Proposition (\ref{thom}) the vector field
\begin{align*}
&w_{ji}(x,y,t)=\dfrac{\partial}{\partial t}+\sum_{j=1}^{n+1}\sum_{i\neq j}
\dfrac{L_{ji}(x,t)}{{{\cal{N}}_{\cal{H}}M_t}}y_{ji}.
\end{align*}
is integrable, which implies ${\cal{C}}^0-{\cal{H}}$-triviality
of $M$.

To prove item $ii)$ let $M(x,t)=M_0(x)+t\Theta(x)$ be a
deformation of $M_0$ with $fil(\theta_{ij})= d_{\max}$,
$\forall\, i,j$ and $t$ is small enough. Analogously to the proof of
item $(i)$, let $p(x,y,t)W(x,y,t)$ be the vector field  where
$W(x,y,t)$ is defined in the Proposition (\ref{thom}) and
$p:\matC^r\times \matC^{(n+1)^2}\times
\matC\longrightarrow\matC$ is a conic bump function ( see
\cite{saia}, Lemma $4$), such that the restriction to $\matC^r\times
{\matC^{(n+1)^2}\times \matC} \setminus \{(0,0,t)\}$ of $p$ 
is smooth and
$$\left\{\begin{array}{l}
  p(x,y,t)=1, \mbox{ for all }(x,y,t)\in\overline{U}\\
  \left.p(x,y,t)=0 \mbox{ in the complement of } V\right.\\
  \left.0\leq |p(x,y,t)|\leq 1 \mbox{ in } V-\overline{U}\right.\\
  p(0,0,t)=0, \mbox{ for all } t
\end{array}\right.$$
where $V$ and $U$ are neighborhoods of the open set
$\{(x,y)\in\mathbb{C}^r\times\mathbb{C}|\,|y|<c\rho(x)\}$ on ${\matC^r\times
\matC^{(n+1)^2}\times \matC\backslash \{(0,0,t)\}}$, where $\rho(x)=[N_{\cal{H}}M_t(x)]^{\frac{1}{2k_1}}$
$$V=\{(x,y,t)\mbox{ such that } |y|\leq c_1\rho(x)\}$$
and $U$ is chosen such that $U\subset\overline{U}\subset V$. Then,
$$|p(x,y,t)W_{ij}(x,y,t)|=\left|\dfrac{L_{ij}(x,t)}{{{\cal{N}}_{\cal{H}}M_t}}\right||py|\leq\left|\dfrac{L_{ij}(x,t)}{({{\cal{N}}_{\cal{H}}M_t)}}\right|{{\cal{N}}_{\cal{H}}M_t}^{1/2k_1},$$
which implies the integrability of the field, hence the
${\cal{C}}^0$-${\cal{H}}$- triviality.

$\hfill \blacksquare$

\vspace{.3cm}

\begin{proposition}
Let $M_0\in Mat_{(n,n+1)}({\cal{O}}_r)$ be a
weighted homogeneous matrix of type $(D,a)\in
Mat_{(n,n+1)}(\matN)\times \matN^r $,
satisfying the condition ${\cal{N}}_{\cal{R}}(M_0(x))\geq
c|x|^{\alpha}$ for constants $c$ and $\alpha$. 
Then for sufficiently small $x$ and $t$, 
deformations $M(x,t)=M_0(x)+t\Theta(x)$, with $fil(\theta_{ij})\geq d_{\max}+1$,
$\forall\, i,j$,  and $d_{max}=\displaystyle{\max_{ij}}\{d_{ij}\}$, are
${\cal{C}}^{0}-{\cal{R}}-$trivial.
\end{proposition}
\newpage
\noindent \textbf{Proof.}
 
Using the proof of the Proposition (\ref{gen}), for each $k,\,l$ fixed we have 
\begin{align}\label{indif}
\Delta^{(q,s)} E_{kl}= (-1)^{l-q}\dfrac{\partial f_l}{\partial x_{\gamma}}G_{sq}^{\gamma}+\dfrac{\partial f_s}{\partial x_{\gamma}}G_{ql}^{\nu}-\dfrac{\partial f_q}{\partial x_{\nu}}G_{sl}^{\nu},
\end{align}
if $q,\,s\neq l$ and
\begin{align}\label{indig}
\Delta^{(l,s)} E_{kl}= \dfrac{\partial f_l}{\partial x_{\gamma}}G_{sl}^{\nu}-\dfrac{\partial f_l}{\partial
x_{\nu}}G_{sl}^{\gamma}
\end{align}
for $j\neq t$, $1\leq k\leq n$, $1\leq j,\,t,\,l\leq n+1$ and $1\leq \gamma,\,\nu\leq r$.

Now multiplying the equation (\ref{indif}) by $|\Delta^{(q,s)}|^{2(\alpha_{qs}-1)}\overline{\Delta}^{(q,s)}$ and adding in $q\neq s$, multiplying the equation (\ref{indig}) by $|\Delta^{(l,s)}|^{2(\alpha_{ls}-1)}\overline{\Delta}^{(l,s)}$ and adding in $s$, we get
\begin{align*}
\sum_{q\neq s}|\Delta^{(q,s)}|^{2\alpha_{qs}}E_{kl}= &\sum_{q>s}|\Delta^{(q,s)}|^{2(\alpha_{qs}-1)}\overline{\Delta}^{(q,s)}\left((-1)^{l-q}\dfrac{\partial f_l}{\partial x_{\gamma}}G_{sq}^{\gamma}+\dfrac{\partial f_s}{\partial x_{\gamma}}G_{ql}^{\nu}-\dfrac{\partial f_q}{\partial x_{\nu}}G_{sl}^{\nu}\right),\\
\sum_{s\neq l}|\Delta^{(l,s)}|^{2\alpha_{ls}}E_{kl}= &\sum_{s\neq l}|\Delta^{(l,s)}|^{2(\alpha_{ls}-1)}\overline{\Delta}^{(l,s)}\left(\dfrac{\partial f_l}{\partial x_{\gamma}}G_{sl}^{\nu}-\dfrac{\partial f_l}{\partial
x_{\nu}}G_{sl}^{\gamma}\right).
\end{align*}

Then it follows that
\begin{align*}
N_{\cal{R}}M\,E_{kl}&=\sum_{\nu=1}^{r}\sum_{\gamma>\nu}^{r}\left[\sum_{s=1}^{n+1}\sum_{q>s}^{n+1}\rho_{qs}\left((-1)^{l-q}\dfrac{\partial f_l}{\partial x_{\gamma}}G_{sq}^{\gamma}+\dfrac{\partial f_s}{\partial x_{\gamma}}G_{ql}^{\nu}-\dfrac{\partial f_q}{\partial x_{\nu}}G_{sl}^{\nu}\right)\right.\left.+\sum_{s\neq l}\rho_{ls}\left(\dfrac{\partial f_l}{\partial x_{\gamma}}G_{sl}^{\nu}-\dfrac{\partial f_l}{\partial x_{\nu}}G_{sl}^{\gamma}\right)\right]
\end{align*}
where $\rho_{qs}=|\Delta^{(q,s)}|^{2(\alpha_{qs}-1)}\overline{\Delta}^{(q,s)}$ and $\rho_{ls}=|\Delta^{(l,s)}|^{2(\alpha_{ls}-1)}\overline{\Delta}^{(l,s)}$.

%
\vspace{.5cm}
Using the Lemma \ref{l4.1} and the matrices above, we have  
\begin{align}
N_{\cal{R}}M\dfrac{\partial M}{\partial t}=\sum_{\nu=1}^{r}\xi_{\nu}\dfrac{\partial M}{\partial x_{\nu}}+\sum_{i=1}^{n+1}B_{ki}C_{li}+\sum_{u\neq k}^{n}S_{ul}R_{ku}\label{control}
\end{align} 
where
\begin{align*}
\xi_{\nu}&=\sum_{\gamma\neq \nu}^r\sum_{k,l}\dfrac{\partial m_{kl}}{\partial t}\left[\sum_{s=1}^{n+1}\sum_{q>s}^{n+1}\rho_{qs}\left(\dfrac{\partial f_s}{\partial x_{\gamma}}cof^q(m_{kl})-\dfrac{\partial f_l}{\partial x_{\gamma}}cof^s(m_{kl})\right)+\sum_{s\neq l}\rho_{ls}\dfrac{\partial f_l}{\partial x_{\gamma}}cof^s(m_{kl})\right]\\
B_{ki}&=\sum_{\nu=1}^{r}\sum_{\gamma>\nu}^{r}\left[\sum_{k,l}\dfrac{\partial m_{kl}}{\partial t}\left(\sum_{s\neq i}^{n+1}\sum_{{\color{red}q\neq i, s}}^{n+1}\rho_{qs}\left(\dfrac{\partial f_s}{\partial x_{\gamma}}\dfrac{\partial cof^s(m_{ki})}{\partial x_{\gamma}}-\dfrac{\partial f_q}{\partial x_{\nu}}\dfrac{\partial cof^s(m_{ki})}{\partial x_{\nu}}\right)+\right.\right.\\&\left.\left.+\sum_{s\neq l,i}\rho_{ls}\left(\dfrac{\partial f_l}{\partial x_{\gamma}}\dfrac{\partial cof^s(m_{ki})}{\partial x_{\nu}}-\dfrac{\partial f_l}{\partial x_{\nu}}\dfrac{\partial cof^s(m_{ki})}{\partial x_{\gamma}}\right)\right)\right]\end{align*}

\begin{align*}
S_{ul}&=\sum_{\gamma>\nu}^{r}\sum_{k,l}\dfrac{\partial m_{kl}}{\partial t}\left[\sum_{s>q}^{n+1}\rho_{qs}\left((-1)^{l-q}\dfrac{\partial f_l}{\partial x_{\gamma}}\sum_{i\neq s,q}\frac{\partial m_{ki}}{\partial x_{\gamma}}cof^{s}_{ki}(m_{uq})+ 
\dfrac{\partial f_s}{\partial x_{\gamma}}\sum_{i\neq q,l}\frac{\partial m_{ki}}{\partial x_{\nu}}cof^{q}_{ki}(m_{ul}) \right.\right. \\ &-\left.\left.\dfrac{\partial f_q}{\partial x_{\nu}}\sum_{i\neq s,l}\frac{\partial m_{ki}}{\partial x_{\nu}}cof^{s}_{ki}(m_{ul})\right)+\sum_{s\neq l} \sum_{i\neq s,l}\rho_{ls}\left(\dfrac{\partial f_l}{\partial x_{\gamma}}\dfrac{\partial m_{ki}}{\partial x_{\nu}}-\dfrac{\partial f_l}{\partial x_{\nu}}\dfrac{\partial m_{ki}}{\partial x_{\gamma}}\right)cof^{s}_{ki}(m_{ul})\right],
\end{align*}

Then,
\begin{align}\label{finally1}
&fil\left(\xi_{\nu}\right)\geq 2k_2+1,\\
&fil\left(B_{ki}\right)\geq 2k_2+1,\\
&fil\left(S_{ul}\right)\geq 2k_2+1.\label{finally2}
\end{align}

Finally, by proposition (\ref{thom}) follows the
${\cal{C}}^0$-${\cal{R}}$-triviality.


$\hfill \blacksquare$

\vspace{0.3cm}
\begin{theorem}
Let $M_0\in Mat_{(n,n+1)}({\cal{O}}_r)$ be a
germ of weighted homogeneous matrix of type $(D,a)\in
Mat_{(n,n+1)}(\matN)\times \matN^r $,
satisfying the condition ${\cal{N}}_{\cal{G}}(M_0(x))\geq
c|x|^{\alpha}$ for constants $c$ and $\alpha$. 
\begin{itemize}
\item[a)] Deformations $M(x,t)=M_0(x)+t\Theta(x)$, with 
$fil(\theta_{ij})\geq d_{\max}+1$, $\forall\, i,j$ and
$d_{max}=\displaystyle{\max_{ij}}\{d_{ij}\}$, are ${\cal{C}}^{0}-{\cal{G}}-$trivial.
\item[b)] Deformations
$M(x,t)=M_0(x)+t\Theta(x)$, with $fil(\theta_{ij})\geq d_{\max}$,
$\forall\, i,j$,  and $d_{max}=\displaystyle{\max_{ij}}\{d_{ij}\}$,
are ${\cal{C}}^{0}-{\cal{G}}-$trivial for  $t$ sufficiently small.
\end{itemize}
\end{theorem}

\noindent \textbf{Proof.}
\begin{itemize}
\item[a)] Since the group ${\cal{C}}^{0}-{\cal{G}}$ is the semi-direct product of the groups ${\cal{C}}^{0}-{\cal{R}}$ and ${\cal{C}}^{0}-{\cal{H}}$, the vector fields are defined as in cases $\mathcal{R}$ and $\mathcal{H}$, and the control function $N_{\mathcal{G}}M$ is defined by $N_{\mathcal{G}}M=N_{\mathcal{R}}^{c_1}M+N_{H}^{c_2}M$ where $c_1$ and $c_2$ are constants such that $N_{\mathcal{G}}M$ is weighted homogeneous.
\item[b)] By Theorem \ref{equi}, if $M_0$ is ${\cal{G}}$-finitely determined,
then ${\cal{N}}_{\cal{G}}(M_0(x))$ satisfies the condition of the
item \textbf{a)}, and we get the result.
\end{itemize}

$\hfill \blacksquare$

\vspace{0.3cm}

%
%
%
%

\begin{example}
Let $M_0=\left(
                           \begin{array}{ccc}
                             z & y & x^3 \\
                             x^2 & z & y \\
                           \end{array}
                         \right)$ weight homogeneous with weights $(3,8,7)$.
Since ${\cal{M}}^3Mat_{n,n+1}({\matC})\subset T{\cal{G}}M_0$ we have
that $M_0$ is ${\cal{G}}$-finitely determined. Therefore, all
deformations $M$ of $M_0$ with filtration of degree greater than
$8$ are ${\cal{C}}^0$-${\cal{G}}$-trivial. Moreover, if $t$ sufficiently small deformations of degree exactly $8$ also are ${\cal{C}}^0$-${\cal{G}}$-trivial.
\end{example}


\begin{thebibliography}{99}
\renewcommand{\refname}{ }

\bibitem{arn}
{\sc V. I. Arnold}, {\it Matrices depending on parameters}, {\sl
Russian Math. Surveys } {\bf 26} , no. 2, 29-46, (1971).

\bibitem{FK} {\sc A. Fr\"uhbis-Kr\"uger}, {\em Classification of Simple Space Curves
Singularities,} {\sl Comm. in Alg.,} {\bf 27 (8)}, (1999),
3993-4013.

\bibitem{FKT} {\sc A. Fr\"uhbis-Kr\"uger}, {\em Moduli Spaces for Space Curve Singularities,}
{\sl Universität Kaiserslautern, Ph. D thesis}, (2000).

\bibitem{FKn} {\sc A. Fr\"uhbis-Kr\"uger, A. Neumer}, {\em Simple Cohen-Macaulay Codimension 2
Singularities,} {\sl Comm. in Alg.}, {\bf 38}, no. 2, 454-495,
(2010).

\bibitem{wall}
{\sc C. T. Wall}, {\it Finite determinacy of smooth map-germs}, {\sl
Bull. London Math. Soc.} {\bf 13} (1981), no 6, 481-539. MR {\bf
83i}:58020.

\bibitem{has}
{\sc G. J. Haslinger}, {\it Families of skew-simmetric matrices},
{\sl University of Liverpool Thesis}, 2001.

\bibitem{Ea1} {\sc J. A. Eagon}, {\em Ideals generated by
subdeterminants of a matrix}, {\sl Thesis}, {\em University of
Chicago}, (1961).


\bibitem{Damon}
{\sc  J. Damon}, {\it The unfolding and determinancy theorems for
subgoups of ${\cal{A}}$ and ${\cal K}$}, {\sl Memoirs of the
American Mathematical Society}, Providence RI, 1984.


\bibitem{Bruce}
{\sc J. W. Bruce}, {\it Families of symmetric matrices}, {\sl Moscow
Math. J.}, {\bf 3}, no 2, (2003), 335-360.

\bibitem{Farid}
{\sc J. W. Bruce and F. Tari}, {\it On families of square matrices},
{\sl Proc. London Math. Soc.} (3), 89, (2004),738-762.


\bibitem{saia}
{\sc M. A. S. Ruas, M. J. Saia}, {\it ${\cal{C}}^l$- determinancy of
weighted homogeneous germs}, {\sl Hokkaido Mathematical Journal}
\textbf{26}, (1997), 89-99.

\bibitem{ruas} {\sc M. A. S. Ruas e J. N. Tomazella}, {\it An infinitesimal criterion for topological triviality of families of sections of analytic varieties},
{\sl Advanced Studies in Pure Mathematcs}, {\sl Proceedings of 12th
MSJ-IRI Symposium "Singularity theory and its applications"}, {\bf
43}, (2006), 421-436.

\bibitem{kuo} {\sc T. C. Kuo}, {\em On ${\cal{C}}^0$-sufficiency of jets of potential functions},
{\sl Topology}, \textbf{8},(1969), 167-171.

\bibitem{Gaft} {\sc T. Gaffney}, {\em Properties of finitely determined germs},
{\sl Brandeis University, Ph. D. thesis}, \textbf{107},(1976).

\bibitem{Bruns} {\sc W. Bruns e U. Vetter}, {\em Determinantal
Rings}, {\sl Springer- Verlang}, New York, (1998).

\bibitem{Bruns1} {\sc W. Bruns e J. Herzog}, {\em Cohen-Macaulay Rings, Revised edition},
{\sl Cambridge University Press}, New York, (1998).

\bibitem{tese} {\sc M. S. Pereira}, {\it Variedades Determinantais e Singularidades de Matrizes},
{\sl Tese de Doutorado}, {\sl ICMC- USP, (2010), Available: \it{http://www.teses.usp.br/teses/disponiveis/55/55135/tde-22062010-133339/pt-br.php}.}

\bibitem{EG} {\sc W. Ebeling and S. M. Gusein-Zade}, \it{On indices of $1$-forms on determinantal singularities},
{\sl Tr. Mat. Inst. Steklova}, {\bf 267}, {\sl (2009), pp. 119-
131}.

\end{thebibliography}
\end{document}